\documentclass[12pt]{amsart}
\usepackage[]{amscd,amsfonts,amsmath,amsxtra,amssymb}
\usepackage{graphics}
\usepackage[all]{xy}
\usepackage{hyperref}
\usepackage[french]{babel}
\usepackage[T1]{fontenc}
\newtheorem{sub}{}[section]
\newtheorem{subsub}{}[sub]
\def\ov#1{\overline{#1}}

\def\coker{\mathop{\rm coker}\nolimits}
\def\Hom{\mathop{\rm Hom}\nolimits}
\def\HHom{\mathop{\mathcal Hom}\nolimits}

\def\Ext{\mathop{\rm Ext}\nolimits}
\def\EExt{\mathop{\mathcal Ext}\nolimits}

\def\Hilb{\mathop{\rm Hilb}\nolimits}
\def\Pic{\mathop{\rm Pic}\nolimits}

\def\imm{\mathop{\rm im}\nolimits}

\def\spec{\mathop{\rm spec}\nolimits}
\def\lra{\longrightarrow}

\def\sigg{\mathop{\hbox{$\displaystyle\sum$}}\limits}

\def\hfl#1#2{\smash{\mathop{\ \hbox to 12mm{\rightarrowfill}}
\limits^{\scriptstyle#1}_{\scriptstyle#2} \ }}
\def\hflb#1#2{\smash{\mathop{\hbox to 12mm{\leftarrowfill}}
\limits^{\scriptstyle#1}_{\scriptstyle#2}}}

\def\m#1{{\hbox{$#1$}}}

\def\ot{\otimes}
\def\og{\leavevmode\raise.3ex\hbox{$\scriptscriptstyle\langle\!\langle$}}
\def\fg{\leavevmode\raise.3ex\hbox{$\scriptscriptstyle\,\rangle\!\rangle$}}

\def\nsp{\lbrace 0\rbrace}

\def\Ssect#1#2{\pagebreak[3]\begin{sub}\label{#2}{\sc\small\small
#1}\rm\medskip}

\def\sepsec{\vskip 1.4cm}
\def\sepsub{\vskip 0.7cm}

\def\sepprop{\vskip 0.5cm}

\def\xmat#1{\[\xymatrix{#1}\]}
\def\flinc{\ar@{^{(}->}}
\def\fleq{\ar@{=}}
\def\flon{\ar@{->>}}
\def\fmaps{\ar@{|-{>}}}
\def\Nligne{\hfil\break}

\newcommand{\C}{{\mathbb C}}

\renewcommand{\P}{{\mathbb P}}

\newcommand{\E}{{\mathbb E}}

\newcommand{\U}{{\mathbb U}}

\renewcommand{\L}{{\mathbb L}}

\newcommand{\kc}{{\mathcal C}}
\newcommand{\kd}{{\mathcal D}}
\newcommand{\ke}{{\mathcal E}}
\newcommand{\kf}{{\mathcal F}}

\newcommand{\ki}{{\mathcal I}}

\newcommand{\kl}{{\mathcal L}}

\newcommand{\ko}{{\mathcal O}}

\newcommand{\ku}{{\mathcal U}}

\begin{document}

\def\refname{R\'ef\'erences}
\def\contentsname{Sommaire}
\def\proofname{D\'emonstration}
\def\abstractname{R\'esum\'e}

\author{Jean--Marc Dr\'{e}zet}
\address{
Institut de Math\'ematiques de Jussieu,
Case 247,
4 place Jussieu,
F-75252 Paris, France}
\email{drezet@math.jussieu.fr}
\title[{Courbes multiples primitives}] {Courbes multiples primitives et 
d\'eformations de courbes lisses}

\begin{abstract}
A {\em primitive multiple curve} is a Cohen-Macaulay scheme $Y$ over $\C$ such 
that the reduced scheme \m{C=Y_{red}} is a smooth curve, and that $Y$ can be 
locally embedded in a smooth surface. In general such a curve $Y$ cannot be 
embedded in a smooth surface. If $Y$ is a primitive multiple curve of 
multiplicity $n$, then there is a canonical filtration \
\m{C=C_1\subset\cdots\subset C_n=Y} \
such that \m{C_i} is a primitive multiple curve of multiplicity $i$. The ideal 
sheaf \m{\ki_C} of $C$ in $Y$ is a line bundle on \m{C_{n-1}}.

Let $T$ be a smooth curve and \m{t_0\in T} a closed point. Let \m{\kd\to T} be 
a flat family of projective smooth irreducible curves, and \m{C=\kd_{t_0}}. 
Then  the $n$-th infinitesimal neighbourhood of $C$ in $\kd$ is a primitive 
multiple curve \m{C_n} of multiplicity $n$, embedded in the smooth surface 
$\kd$, and in this case \m{\ki_C} is the 
trivial line bundle on \m{C_{n-1}}. Conversely, we prove that every projective 
primitive multiple curve \m{Y=C_n} such that \m{\ki_C} is the trivial line 
bundle on \m{C_{n-1}} can be obtained in this way.
\end{abstract}

\maketitle
\tableofcontents

{\em Mathematics Subject Classification :} 14H10, 14C05

\vskip 1cm

\section{Introduction}

\Ssect{Courbes multiples primitives}{cmpr}

Une {\em courbe primitive} est une vari\'et\'e alg\'ebrique complexe
$Y$, de Cohen-Macaulay et telle que la sous-vari\'et\'e r\'eduite associ\'ee 
\m{C=Y_{red}} soit une courbe lisse irr\'eductible, et que tout point ferm\'e 
de $Y$ poss\`ede un voisinage pouvant \^etre plong\'e dans une surface lisse. 
Ces courbes ont \'et\'e d\'efinies et \'etudi\'ees par C.~B\u anic\u a et
O.~Forster dans \cite{ba_fo}. On s'int\'eresse plus particuli\`erement ici au
cas o\`u $C$, et donc $Y$, sont projectives.

Soient $P$ un point ferm\'e de $Y$, et $U$ un voisinage de $P$ pouvant
\^etre plong\'e dans une surface lisse $S$. Soit $z$ un \'el\'ement de
l'id\'eal maximal de l'anneau local \m{\ko_{S,P}} de $S$ en $P$ engendrant
l'id\'eal de $C$ dans cet anneau. Il existe alors un unique entier $n$,
ind\'ependant de $P$, tel que l'id\'eal de $Y$ dans \m{\ko_{S,P}} soit
engendr\'e par \m{(z^n)}. Cet entier $n$ s'appelle la {\em multiplicit\'e} de
$Y$. Si \m{n=2} on dit que $Y$ est une {\em courbe double}.
Il existe une filtration canonique
\[C=C_1\subset\cdots\subset C_n=Y\ ,\]
o\`u au voisinage de chaque point $P$ l'id\'eal de \m{C_i} dans \m{\ko_{S,P}}
est \m{(z^i)}. Donc \m{C_i} est une courbe multiple primitive de 
multiplicit\'e $i$.

Soit \m{\ki_C} le faisceau d'id\'eaux de $C$ dans $Y$. Alors le faisceau
conormal de $C$, \m{L=\ki_C/\ki_C^2} est un fibr\'e en droites sur $C$, dit
{\em associ\'e} \`a $Y$, et \m{\ki_C} est un fibr\'e en droites sur 
\m{C_{n-1}}.

Les exemples les plus simples de courbes multiples primitives sont les courbes 
de Cohen-Macaulay de courbe r\'eduite associ\'ee lisse et plong\'ees dans une 
surface lisse. En particulier, soit \m{L\in\Pic(C)}. Alors le $n$-i\`eme 
voisinage infinit\'esimal de $C$ dans \m{L^*} est une courbe primitive de 
multiplicit\'e $n$ et de fibr\'e en droites associ\'e $L$. On l'appelle la 
{\em courbe primitive triviale} de fibr\'e associ\'e $L$.

\end{sub}

\sepsub

\Ssect{Courbes d\'efinies par des familles de courbes lisses}{cdfl}

En g\'en\'eral une courbe multiple primitive ne peut pas
\^etre plong\'ee dans une surface lisse. D'apr\`es \cite{ba_ei},
theorem 7.1, la seule courbe double non triviale de courbe r\'eduite
associ\'ee \m{\P_1} pouvant \^etre plong\'ee dans une surface lisse est la
courbe double d\'eduite d'une conique plane. Or il existe beaucoup d'autres 
courbes doubles de courbe r\'eduite associ\'ee \m{\P_1}.

On s'int\'eresse dans cet article aux courbes multiples d\'efinies par les 
familles de courbes lisses. Ce sont des cas particuliers de courbes multiples 
primitives plong\'ees dans une surface lisse. Soient $T$ une courbe lisse (ou 
un germe de courbe lisse), \m{t_0} un point ferm\'e de $T$, et \m{\kd\to T} 
une famille plate de courbes projectives lisses et irr\'eductibles 
param\'etr\'ee par $T$. Soit \ \m{C=\kd_{t_0}}. Pour tout entier \m{n\geq 2}, 
le n-i\`eme voisinage infinit\'esimal \m{C_n} de $C$ dans $\kd$ est une courbe 
primitive de multiplicit\'e $n$, de courbe lisse associ\'ee $C$. Pour cette 
courbe primitive le faisceau d'id\'eaux \m{\ki_C} est le fibr\'e trivial sur 
\m{C_{n-1}}. On dit d'une telle courbe qu'elle {\em provient d'une famille de 
courbes lisses}.

Le principal r\'esultat de cet article est le

\sepprop

\begin{subsub}\label{theo}{\bf Th\'eor\`eme : } Soit \m{C_n} une courbe 
multiple primitive de multiplicit\'e $n$, de courbe lisse sous-jacente $C$. 
Alors \m{C_n} provient d'une famille de courbes lisses si et seulement si le 
faisceau d'id\'eaux de $C$ dans \m{C_n} est trivial sur \m{C_{n-1}}.
\end{subsub}

\sepprop

La d\'emonstration repose sur la param\'etrisation des courbes multiples 
primitives. Les courbes doubles ont \'et\'e d\'ecrites dans \cite{ba_ei}. Le 
th\'eor\`eme \ref{theo} pour les courbes doubles en d\'ecoule ais\'ement.
On fait ensuite une d\'emonstration par r\'ecurrence sur $n$, en s'appuyant 
sur la description des courbes de multiplicit\'e \m{n>2} donn\'ee dans 
\cite{dr1}.

Les courbes multiples primitives provenant de familles de courbes lisses 
peuvent aussi \^etre construites de la fa\c con suivante : on part d'une 
famille plate \m{\pi:\kd\to S} de courbes lisses param\'etr\'ee par une 
vari\'et\'e lisse $S$. Soit \m{s_0\in S} un point ferm\'e tel que 
\m{\kd_{s_0}\simeq C}. Soit \ \m{Z_n=\spec(\C[t]/(t^n))\hookrightarrow S} \ un 
plongement tel que l'image du point ferm\'e de \m{Z_n} soit \m{s_0}. Alors 
\m{\pi^{-1}(Z_n)} est une courbe multiple primitive de multiplicit\'e $n$ 
provenant d'une famille de courbes lisses (cela revient \`a inclure \m{Z_n} 
dans une courbe lisse de $S$).

Soit \m{C_n} une courbe multiple primitive de multiplicit\'e $n$ telle que le 
faisceau d'id\'eaux \m{\ki_C} soit trivial sur \m{C_{n-1}}. On a une suite 
exacte canonique de faisceaux coh\'erents sur \m{C_{n-1}}
\[0\lra\ko_C\lra\Omega_{C_n\mid C_{n-1}}\lra\Omega_{C_{n-1}}\lra 0 .\]
A tout prolongement de \m{C_{n-1}} en courbe de multiplicit\'e $n$ 
correspond un certain type d'extension de \m{\Omega_{C_{n-1}}} par \m{\ko_C} 
sur \m{C_{n-1}} dont le terme du milieu est localement libre. On montre que 
cette correspondance est bijective (un r\'esultat analogue a \'et\'e obtenu 
dans \cite{dr1}, 6-).

On consid\`ere ensuite une famille plate de courbes lisses \m{\kd\to S} comme 
pr\'ec\'edemment, telle que le morphisme de Koda\"\i ra-Spencer \ 
\m{KS:T_{s_0}S\to H^1(C,T_C)} \ soit surjectif. On suppose que \m{C_{n-1}} 
provient d'un plongement \m{Z_{n-1}\hookrightarrow S}. On fait ensuite le lien 
entre les prolongements de ce plongement \`a \m{Z_n\hookrightarrow S} et les 
extensions de \m{\Omega_{C_{n-1}}} par \m{\ko_C}.

Soit $g$ le genre de $C$, et supposons que \m{g\geq 2}. Alors les 
prolongements de \m{C_{n-1}} en courbes de multiplicit\'e $n$ telles que 
\m{\ki_C} soit trivial sur \m{C_{n-1}} forment un espace de dimension \m{3g-3}
si \m{C_2} est triviale, et \m{3g-4} sinon.

\end{sub}

\sepsub

\Ssect{Motivation}{moti}

Soient $S$ une courbe lisse, \m{\pi:\kc\to S} une famille plate de courbes
projectives, et \m{s_0\in S} tel que \m{\kc_{s_0}} soit une courbe multiple
primitive de multiplicit\'e \m{n>0}, de courbe lisse associ\'ee $C$. Soient
\m{\ko(1)} un fibr\'e en droites tr\`es ample sur $\kc$ et $P$ un polyn\^ome en
une variable \`a coefficients rationnels. On a alors une {\em vari\'et\'e de
modules relative}
\[\rho:{\bf M}_{\ko(1)}(P)\lra S \]
des faisceaux semi-stables sur les fibres de $\pi$ de polyn\^ome de Hilbert $P$
(cf. \cite{si}). En g\'en\'eral $\rho$ n'est pas plat (cf. \cite{in2}). Par
exemple, si les courbes \m{\kc_s}, \m{s\not=s_0}, sont lisses, \m{\ko_C}, vu
comme faisceau sur \m{\kc_{s_0}}, est stable, mais ne se d\'eforme pas en
faisceaux stables sur les autres fibres \m{\kc_s}.

Je conjecture que si les fibres \m{\kc_s}, \m{s\not=s_0}, ont exactement $n$
composantes irr\'eductibles qui sont lisses, alors $\rho$ est plat. Les
faisceaux (semi-)stables sur ce genre de courbe sont bien connus (cf.
\cite{tei2}). Il reste cependant \`a classifier les courbes multiples primitives
qui sont des d\'eformations de courbes r\'eductibles \`a composantes lisses.

Il est facile de voir que les courbes multiples provenant de d\'eformations de 
courbes lisses \'etudi\'ees ici sont des cas particuliers de d\'eformations de 
courbes r\'eductibles \`a composantes lisses. C'est l'exemple le plus simple, et 
la classification est donc donn\'ee dans ce cas.

\end{sub}

\newpage

\section{Pr\'eliminaires}\label{prelim}

\Ssect{Cohomologie de \v Cech}{cech}

On donne ici plusieurs fa\c cons de repr\'esenter des faisceaux ou des classes
de cohomologie en utilisant diff\'erentes variantes de la cohomologie de \v 
Cech.

\sepprop

\begin{subsub}\label{fai_rec} D\'efinition de faisceaux par recollements -- \rm
Soient $X$ une vari\'et\'e alg\'ebrique et \m{(X_i)_{i\in I}} un revouvrement
ouvert de $X$. Pour tout \m{i\in I}, soient \m{U_i} une vari\'et\'e
alg\'ebrique et \ \m{\alpha_i:X_i\to U_i} \ isomorphisme. Pour tous \m{i,j\in 
I} distincts, soit \ \m{U_{ij}^{(i)}=\alpha_i(X_{ij})}.

La vari\'et\'e $X$ est obtenue en ``recollant'' les vari\'et\'es \m{U_i} au
moyen des isomorphismes
\[\alpha_{ij}=\alpha_j\circ\alpha_i^{-1}:U_{ij}^{(i)}\lra U_{ij}^{(j)} .\]

Soit $\kf$ un faisceau coh\'erent sur $X$. On en d\'eduit \ 
\m{\kf_i=(\alpha_i^{-1})^*(\kf)}, faisceau coh\'erent sur \m{U_i}. On a des
isomorphismes canoniques
\xmat{\Theta_{ij}:\kf_{i\mid U_{ij}^{(i)}}\ar[r]^-\simeq & 
(\alpha_{ij})^*(\kf_{j\mid U_{ij}^{(j)}})}
tels que
\[\Theta_{ik} \ = \ (\alpha_{ij})^*(\Theta_{jk})\circ\Theta_{ij}\]
sur \ \m{U_{ijk}^{(i)}=\alpha_i(X_{ijk})}.

R\'eciproquement, \'etant donn\'es des faisceaux \m{\kf_i} sur \m{U_i} et des
isomorphismes \m{\Theta_{ij}} poss\'edant les propri\'et\'es pr\'ec\'edentes,
on construit ais\'ement un faisceau $\kf$ sur $X$ en recollant les
\m{\alpha_i^*(\kf_i)}.
\end{subsub}

\sepprop

\begin{subsub}\label{coho_cl} D\'efinition de classes de cohomologie -- \rm
On conserve les notations de \ref{fai_rec}. Soit \m{u\in H^1(X,\kf)},
repr\'esent\'e par un cocycle \m{(u_{ij})}, avec \m{u_{ij}\in
H^0(X_{ij},\kf)}. Soit \m{v_{ij}\in H^0(U_{ij}^{(i)},\kf_i)} \ correspondant
\`a \m{u_{ij}}. Alors on a dans \m{H^0(U_{ijk}^{(i)},\kf_i)}
\[v_{ik} \ = \ v_{ij} + \Theta_{ij}^{-1}(v_{jk}) . \]
R\'eciproquement, toute famille \m{(v_{ij})} v\'erifiant les \'egalit\'es
pr\'ec\'edentes d\'efinit un \'el\'ement de \m{H^1(X,\kf)}.

Soit $\kl$ un fibr\'e en droites sur $X$, repr\'esent\'e par une famille
\m{(\lambda_{ij})}, \m{\lambda_{ij}\in H^0(U_{ij}^{(i)},\ko_{U_i}^*)}, comme
dans \ref{fai_rec} (le faisceau sur \m{U_i} pour $\kl$ \'etant donc
\m{\ko_{U_i}}). On a donc des ismorphismes canoniques \ \m{
(\alpha_i^{-1})^*(\kl)\simeq\ko_{U_i}}.

Soit \m{\mathbf{u}\in H^1(X,\kf\ot\kl)}, repr\'esent\'e par un
cocycle \m{(\mathbf{u}_{ij})}, avec \ \m{\mathbf{u}_{ij}\in
H^0(X_{ij},\kf\ot\kl)}. On d\'eduit de \m{\mathbf{u}_{ij}} et des isomorphismes
\m{\kf_i\simeq(\alpha_i^{-1})^*(\kf)},
\m{\ko_{U_j}\simeq(\alpha_j^{-1})^*(\kl)} (bien noter que $i$ est utilis\'e
pour $\kf$ et $j$ pour $\kl$) un \'el\'ement \m{\mathbf{v}_{ij}} de
\m{H^0(U_{ij}^{(i)},\kf_i)}. Plus pr\'ecis\'ement, \m{\mathbf{v}_{ij}} est
obtenu au moyen d'un isomorphisme \ \m{(\alpha_i^{-1})^*(\kf\ot\kl)\simeq\kf_i}
\ sur \m{U_{ij}^{(i)}} : si \m{x\in X_{ij}}, \m{f\in\kf_x}, \m{\ell\in\kl_x},
on en d\'eduit \m{y=\alpha_i(x)\in U_{ij}^{(i)}}, \m{f_i\in\kf_{i,y}},
\m{\ell_j\in\ko_{U_j,\alpha_j(x)}}, \m{\ell_i=\Theta_{ij}^*(\ell_j)\in
\ko_{U_j,\alpha_j(x)}}. D'o\`u \m{\ell_if_i\in\kf_{i,y}}.

On a alors dans \m{H^0(U_{ijk}^{(i)},\kf_i)}
\[\mathbf{v}_{ik} \ = \ \Theta_{ij}^*(\lambda_{jk})\mathbf{v}_{ij} +
\Theta_{ij}^{-1}(\mathbf{v}_{jk}) .   \]
R\'eciproquement, toute famille \m{(\mathbf{v}_{ij})} v\'erifiant les
\'egalit\'es pr\'ec\'edentes d\'efinit un \'el\'ement de \m{H^1(X,\kf\ot\kl)}.
\end{subsub}

\end{sub}

\sepsub

\Ssect{Courbes multiples primitives}{def_nota}

Soit $C$ une courbe lisse irr\'eductible et projective.
Soit $Y$ une courbe multiple primitive de multiplicit\'e $n$, de courbe lisse 
associ\'ee $C$ et de fibr\'e en droites sur $C$ associ\'e $L$. Soit
\[C=C_1\subset\cdots\subset C_n=Y\]
la filtration canonique. Pour \m{2\leq i\leq n}, \m{C_i} est une courbe 
multiple primitive de multiplicit\'e $i$. On notera \ \m{\ko_i=\ko_{C_i}}.

\sepprop

\begin{subsub} Structure locale -- \rm
D'apr\`es \cite{dr1}, th\'eor\`eme 5.2.1, l'anneau \m{\ko_{Y,P}} est isomorphe 
\`a \Nligne \m{\ko_{C,P}\ot(\C[t]/(t^n))}.
\end{subsub}

\sepprop

\begin{subsub}\label{mult_2} Courbes de multiplicit\'e 2 - \rm (cf. 
\cite{ba_ei}) Deux courbes multiples primitives \m{C_2}, \m{C'_2} de courbe 
lisse associ\'ee $C$ sont dites {\em isomorphes} s'il existe un isomorphisme 
\m{C_2\simeq C'_2} induisant l'identit\'e de $C$. Dans ce cas les fibr\'es en 
droites sur $C$ associ\'es \`a \m{C_2}, \m{C'_2} sont isomorphes.

Soient $C$ une courbe lisse irr\'eductible et \m{L\in\Pic(C)}.
Soit \m{C_2} une courbe multiple primitive de multiplicit\'e 2, de courbe 
lisse associ\'ee $C$ et de fibr\'e en droites sur $C$ associ\'e $L$.
Soit \ \m{E=(\Omega_{C_2\mid C})^*}, qui est un fibr\'e vectoriel de rang 2 
sur $C$. On a une suite exacte canonique
\[0\lra T_C\lra E\lra L^*\lra 0\]
associ\'ee \`a un \'el\'ement de 
\m{\Ext^1_{\ko_C}(L^*,T_C)=H^1(C,T_C\ot L)}. R\'eciproquement, on montre que 
pour toute suite exacte $(S)$ du type pr\'ec\'edent, il existe une courbe 
\m{C_2} dont la suite exacte associ\'e est $(S)$.

On montre que deux courbes \m{C_2}, \m{C'_2}, de courbe lisse et fibr\'e 
associ\'es $C$ et $L$, sont isomorphes si et seulement si les \'el\'ements 
correspondants de \m{H^1(T_C\ot L)} sont proportionnels. La courbe 
correspondant \`a 0 est la courbe triviale. Les courbes doubles non triviales 
de courbe lisse associ\'ee $C$ et de fibr\'e en droites associ\'e $L$ sont 
donc param\'etr\'ees par l'espace projectif \m{\P(H^1(T_C\ot L))}.
\end{subsub}

\sepprop

\begin{subsub} Filtration canonique d'un faisceau coh\'erent -- \rm
Le faisceau \m{\ki_C} est un fibr\'e en droites sur \m{C_{n-1}}. Il existe
d'apr\`es \cite{dr2}, th\'eor\`eme 3.1.1, un fibr\'e en droites $\L$ sur
\m{C_n} dont la restriction \`a \m{C_{n-1}} est \m{\ki_C}. On a alors, pour
tout faisceau de \m{\ko_n}-modules $\ke$ un morphisme canonique
\begin{equation}\label{equ4}\ke\ot\L\lra\ke\end{equation}
qui en chaque point ferm\'e $P$ de $C$ est la multiplication par $z$.

Soit $\ke$ un faisceau coh\'erent sur \m{C_n}. Pour tout entier \m{i\geq 0} on
note \m{\ke^{(i)}} le noyau du morphisme \ \m{\ke\to\ke\ot\L^{-i}} \ d\'eduit
de (\ref{equ4}). Au point $P$, \m{\ke_P^{(i)}} est donc le sous-module de 
\m{\ke_P} constitu\'e des \'el\'ements annul\'es par \m{z^i}. La filtration
\[\ke^{(0)}=0\subset\ke^{(1)}\subset\cdots\subset\ke^{(n)}=\ke \]
s'appelle la {\em seconde filtration canonique} de $\ke$. Ses gradu\'es sont
des faisceaux sur $C$ (cf. \cite{dr4}, 3.).
\end{subsub}

\end{sub}

\sepsub

\Ssect{Fibr\'e d\'eterminant}{fib_det}

On conserve les notations de \ref{def_nota}.

Soit $\ke$ un faisceau coh\'erent sur \m{C_n}. On pose
\[{\det}_C(\ke) \ = \ \bigotimes_{i=1}^n\det(\ke^{(i)}/\ke^{(i-1)}) \ .\]

\sepprop

\begin{subsub}{\bf Proposition : }\label{pr4} Soit \
\m{\kf_0\subset\kf_1\subset\cdots\subset\kf_m=\ke} \ une filtration de $\ke$
dont les gradu\'es sont des faisceaux sur $C$. Alors on a
\[{\det}_C(\ke) \ = \ \bigotimes_{i=1}^m\det(\kf_{i}/\kf_{i-1}) \ .\]
\end{subsub}

\begin{proof}
On a \m{\kf_1\subset\ke^{(1)}}. Soit \m{F=\ke^{(1)}/\kf_1}.
On a un diagramme commutatif avec lignes et colonnes exactes
\xmat{  &                       &                   & 0\ar[d]\\
        & 0\ar[d]               &                   & F\ar[d] \\
0\ar[r] & \kf_1\ar[r]\ar[d]     & \ke\ar[r]\fleq[d] & \ke/\kf_1\ar[r]\ar[d] & 
0\\
0\ar[r] & \ke^{(1)}\ar[r]\ar[d] & \ke\ar[r]               &
\ke/\ke^{(1)}\ar[r]\ar[d] & 0\\
        & F\ar[d]               &                   & 0\\
        & 0 
}
Supposons que le r\'esultat soit vrai pour le faisceau \m{\ke/\kf_1}. On
consid\`ere la filtration de \m{\ke/\kf_1}
\[0\subset F\subset\ke^{(2)}/\kf_1\cdots\subset\ke/\kf_1 . \]
On en d\'eduit que
\[\bigotimes_{i=2}^m\det(\kf_{i}/\kf_{i-1})={\det}_C(\ke/\kf_1)=\det(F)\ot{\det
}_C(\ke/\ke^{(1)}) .\]
Donc
\begin{eqnarray*}\bigotimes_{i=1}^m\det(\kf_{i}/\kf_{i-1}) & = & 
\det(\kf_1)\ot\det(F)\ot{\det}_C(\ke/\ke^{(1)})\\
& = & \det(\ke^{(1)})\ot{\det}_C(\ke/\ke^{(1)})\\
& = & {\det}_C(\ke) .
\end{eqnarray*}
Donc le r\'esultat est vrai pour $\ke$.

Il suffit donc de prouver le r\'esultat pour \m{\ke/\kf_1}, qui lui-m\^eme
d\'ecoule de celui pour \m{\ke/\kf_2}. On se ram\`ene ainsi \`a prouver le
r\'esultat pour $\kf_m$. Mais c'est bien connu pour les faisceaux sur $C$.
\end{proof}

\sepprop

\begin{subsub}{\bf Corollaire : }\label{pr5} Soit \ 
\m{0\to\ke'\to\ke\to\ke''\to 0} \ une suite exacte de faisceaux coh\'erents
sur \m{C_n}. Alors on a
\[{\det}_C(\ke) \ = \ {\det}_C(\ke')\ot{\det}_C(\ke'') .\]
\end{subsub}

\end{sub}

\sepsub

\Ssect{Construction des courbes multiples primitives}{prim_const}

(cf. \cite{dr1})

Soit \m{n\geq 2} un entier. On pose \ \m{Z_n=\spec({\C[t]/(t^n)})}. Pour toute
vari\'et\'e alg\'ebrique affine \m{U=\spec(A)}, on a \ \m{U\times
Z_n=\spec(A[t]/(t^n))}. Si \m{u\in A[t]/(t^n)}, on notera \m{u_i} le
coefficient de \m{t^i} dans $u$.

Soit \m{C_n} une courbe multiple primitive de multiplicit\'e $n$, de courbe
r\'eduite associ\'ee $C$ et de fibr\'e en droites sur $C$ associ\'e $L$. Soit
\m{(U_i)} un recouvrement ouvert de $C$ tel que chaque $U_i$ soit affine
(c'est-\`a-dire distinct de $C$ si $C$ est projective), et que les
restrictions \m{\omega_{C\mid U_i}} et \m{L_{\mid U_i}} soient triviales.
Alors \m{C_n} peut se construire en recollant les vari\'et\'es \ \m{U_i\times
Z_n} \ au moyen d'automorphismes \m{\sigma_{ij}} des \ \m{U_{ij}\times Z_n} \
laissant \m{U_{ij}} invariant, la famille \m{(\sigma_{ij})} v\'erifiant la
relation de cocycle
\[\sigma_{ik} \ = \ \sigma_{jk}\circ\sigma_{ij} . \]
Soit $U$ un ouvert affine de $C$ tel que \m{\omega_{C\mid U}} soit trivial, et
soit \m{x\in\ko_C(U)} tel que \m{dx} engendre \m{\omega_{C\mid U}}. Alors les
automorphismes d'alg\`ebres de \ \m{\ko(U\times Z_n)\simeq\ko_C(U)[t]/(t^n)} \
sont de la forme \m{\phi_{\mu,\nu}}, avec \m{\mu,\nu\in\ko_C(U)[t]/(t^{n-1})},
$\nu$ inversible, o\`u, pour tout \m{\alpha\in\ko_C(U)},
\[\phi_{\mu,\nu}(\alpha) \ = \ \sigg_{i=0}^{n-1}\frac{1}{i!}(\mu t)^i
\frac{\partial^i\alpha}{\partial^ix},\]
et \ \m{\phi_{\mu,\nu}(t)=\nu t} . Remarquons que \m{\phi_{\mu,\nu}} est
enti\`erement d\'etermin\'e par
\[\phi_{\mu,\nu}(x) \ = \ x+\mu t \quad\quad \text{et} \quad\quad
\phi_{\mu,\nu}(t) \ = \  \nu t \ . \]
En toute rigueur on devrait noter ce morphisme \m{\phi_{\mu,\nu}^{dx}}, car
$\mu$ d\'epend du choix de $x$.

Supposons que \m{\sigma_{ij}^{-1}} corresponde \`a
\m{\phi_{\mu_{ij},\nu_{ij}}}. On a alors les relations de cocycle
\[\phi_{\mu_{ik},\nu_{ik}} \ = \ \phi_{\mu_{jk},\nu_{jk}}\circ
\phi_{\mu_{ij},\nu_{ij}} .\]

La famille \m{(\nu_{ij})} n'est pas en g\'en\'eral un cocycle de
\m{\ko_{n-1}^*}, mais \m{(\nu_{ij,0})} est un cocycle de \m{\ko_C^*}, et
l'\'el\'ement induit de \m{H^1(\ko_C^*)} est le fibr\'e en droites $L$.

\sepprop

\begin{subsub}\label{prop_phi}
Propri\'et\'es des automorphismes \m{\phi_{\mu,\nu}} -- \rm Si
\m{\mu,\mu'\in\ko_C(U)[t]/(t^{n-1})} et\Nligne
\m{\nu,\nu'\in\ko_C(U)[t]/(t^{n-1})} sont inversibles, on a
\[\phi_{\mu'\nu'}\circ\phi_{\mu\nu} \ = \ \phi_{\mu"\nu"} ,\]
avec
\[\mu"=\mu'+\nu'\phi_{\mu',\nu'}(\mu) , \ \ \ \
\nu"=\nu'\phi_{\mu',\nu'}(\nu) .\]
On a \ \m{\phi_{\mu\nu}^{-1}=\phi_{\ov{\mu},\ov{\nu}}}, avec
\[\ov{\mu} \ = \ -\phi_{\mu\nu}^{-1}(\frac{\mu}{\nu}) , \quad\quad
\ov{\nu} \ = \ \phi_{\mu\nu}^{-1}(\frac{1}{\nu}) .\]
\end{subsub}

\sepprop

\begin{subsub}\label{prolong} Prolongement de courbes multiples - \rm On
suppose que \m{n\geq 3}. Soit \m{C_{n-1}} une courbe multiple primitive de
multiplicit\'e \m{n-1}, de courbe r\'eduite associ\'ee $C$ et de fibr\'e en
droites sur $C$ associ\'e $L$. On note comme dans \ref{mult_2} \ 
\m{E=(\Omega_{C_2\mid C})^*}.

Soit \m{C_n} une courbe multiple de multiplicit\'e $n$ dont la courbe de
multiplicit\'e \m{n-1} sous-jacente est \m{C_{n-1}}. On dit qu'une telle 
courbe est un {\em prolongement} de \m{C_{n-1}} en courbe multiple primitive 
de multiplicit\'e $n$. Deux tels prolongements \m{C_n,C'_n} sont dits {\em 
isomorphes} s'il existe un isomorphisme \m{C_n\simeq C'_n} induisant 
l'identit\'e de \m{C_{n-1}}. On d\'efinit dans \cite{dr1}, \m{C_n} \'etant 
donn\'e, une param\'etrisation des classes d'isomorphisme de 
prolongements de \m{C_{n-1}} en courbe multiplicit\'e $n$ par \m{H^1(E\ot 
L^{n-1})} (\m{C_n} correspondant a $0$).

On peut retrouver ce r\'esultat \`a partir des cocycles. Supposons que \m{C_n}
soit obtenue \`a partir d'une famille \m{(\sigma_{ij})} comme
pr\'ec\'edemment. Une autre extension \m{C'_n} de \m{C_{n-1}} provient d'une
autre famille \m{(\sigma'_{ij})}, o\`u \
\m{{\sigma'}_{ij}^{-1}=\phi_{\mu'_{ij},\nu'_{ij}}}, \m{\mu'_{ij}}, 
\m{\mu'_{ij}} \'etant de la forme
\[\mu'_{ij} \ = \ \mu_{ij}+\alpha_{ij}t^{n-2} , \quad\quad
\nu'_{ij} \ = \ \nu_{ij}+\beta_{ij}t^{n-2} ,  \]
avec \ \m{\alpha_{ij},\beta_{ij}\in\ko_C(U_{ij})}.
En utilisant \ref{prop_phi}, on voit que la relation \
\m{\sigma'_{ik}=\sigma'_{jk}\circ\sigma'_{ij}}, compte tenu du fait que \
\m{\sigma_{ik}=\sigma_{jk}\circ\sigma_{ij}}, \'equivaut \`a l'\'egalit\'e
\begin{equation}\label{equ2}
\begin{pmatrix}\alpha_{ik} \\ \beta_{ik}\end{pmatrix} \ =  \
\nu_{ij,0}^{n-1}\begin{pmatrix}\alpha_{jk} \\ \beta_{jk}\end{pmatrix} +
\begin{pmatrix}1 & \mu_{jk,0} \\ 0 & \nu_{jk,0}\end{pmatrix}
\begin{pmatrix}\alpha_{ij} \\ \beta_{ij}\end{pmatrix} \ .
\end{equation}
D'apr\`es \ref{coho_cl} et la construction de \m{\Omega_{C_2\mid C}} donn\'ee
dans \ref{can_coc}, ces relations montrent que \m{((\alpha_{ij},\beta_{ij}))}
d\'efinit un \'el\'ement $u$ de \m{H^1(E\ot L^{n-1})}.
En reprenant la param\'etrisation des extensions de \m{C_n} donn\'ee dans
\cite{dr1} on voit ais\'ement que $u$ est pr\'ecis\'ement l'\'el\'ement de
\m{H^1(E\ot L^{n-1})} correspondant \`a \m{C'_n}.
\end{subsub}

\sepprop

\begin{subsub}\label{nota_1}{\bf Notation : }
On notera \ \m{C'_n=C_n(u)} (et donc \m{C_n=C_n(0)}).
\end{subsub}
\end{sub}

\sepprop

\begin{subsub}\label{prol_cst} Prolongements \`a faisceau d'id\'eaux constant
-- \rm Le faisceau d'id\'eaux \m{\ki_{C,C_n}} de $C$ dans \m{C_n} est un 
fibr\'e en droites sur \m{C_{n-1}}. Il est construit en recollant les id\'eaux 
\m{(t)} des faisceaux \m{\ko_{U_i\times Z_n}} au moyen des isomorphismes 
\m{\phi_{\ov{\mu}_{ij},\ov{\nu}_{ij}}}, o\`u \m{\ov{\mu}_{ij},\ov{\nu}_{ij}} 
sont les images de \m{\mu_{ij},\nu_{ij}} respectivement dans 
\m{\ko_C(U_{ij})[t]/(t^{n-2})}. On en d\'eduit ais\'ement que les extensions 
\m{C_n(u)} de \m{C_{n-1}} telles que \ \m{\ki_{C,C_n(u)}=\ki_{C,C_n}} \ sont
celles correspondant aux $u$ d\'efinis par des familles 
\m{((\alpha_{ij},\beta_{ij}))} o\`u les \m{\beta_{ij}} sont nuls. 
C'est-\`a-dire que ce sont les courbes \m{C_n(u)}, o\`u $u$ appartient \`a 
l'image de \m{H^1(T_C\ot L^{n-1})} dans \m{H^1(E\ot L^{n-1})}.
\end{subsub}

\sepprop

\begin{subsub}\label{ic_tr} Courbes multiples \`a faisceau d'id\'eaux trivial 
-- \rm Si \m{\ki_{C,C_n}} est trivial (sur \m{C_{n-1}}), le fibr\'e en droites 
sur $C$ associ\'e \`a \m{C_n} est \m{\ko_C}. Soit \m{C_2} une courbe double de 
fibr\'e en droites associ\'e \m{\ko_C}. Si \ \m{E=(\Omega_{C_2\mid C})^*}, on 
a une suite exacte
\[0\lra T_C\lra E\lra\ko_C\lra 0 ,\]
correspondant \`a \ \m{\sigma\in H^1(T_C)}. Soit \ \m{D=\C\sigma\subset 
H^1(T_C)}. On a vu dans \ref{mult_2} que \m{C_2} est enti\`erement 
d\'etermin\'ee par $D$.

On suppose maintenant que \m{n\geq 3}. Soit \m{C_n} une courbe multiple 
primitive de multiplicit\'e $n$ de courbe lisse associ\'ee $C$, telle que  
\m{\ki_{C,C_n}} soit trivial sur \m{C_{n-1}}. On a vu pr\'ec\'edemment que les 
prolongements de \m{C_{n-1}} en courbe \m{C'_n} de multiplicit\'e $n$ telle 
que  \m{\ki_{C,C'_n}} soit trivial \'etaient param\'etr\'ees par l'image de 
\m{H^1(T_C)} dans \m{H^1(E)}, c'est-\`a-dire par \m{H^1(T_C)/D} .
\end{subsub}

\sepsub

\Ssect{Un lemme sur les extensions}{ex_mod}

Soient $X$ une vari\'et\'e alg\'ebrique, et $G$, $F$, $\E$, $\ku$ des
faisceaux coh\'erents sur $X$.
On suppose qu'on a une suite exacte
\xmat{0\ar[r] & G\oplus F\ar[r]^-{\beta} & \E\ar[r] & \ku\ar[r] & 0 .}
Pour tout morphisme \ \m{\psi:G\to F}, on note \m{\ke_\psi} le conoyau du
morphisme compos\'e
\xmat{G\ar[rr]^-{(I_G,-\psi)} & & G\oplus F\ar[r]^\beta & \E .}
On a donc un diagramme commutatif avec lignes et colonnes exactes
\xmat{& & & & 0\ar[d] \\ & 0\ar[d] & & & F\ar[d] \\
0\ar[r] & G\ar[rr]^-{\beta(I_G,-\psi)}\ar[d]^{(I_G,-\psi)} & & 
\E\fleq[d]\ar[r] & \ke_\psi\ar[r]\ar[d]& 0\\
0\ar[r] & G\oplus F\ar[rr]^-\beta\ar[d]^{(\psi,I_F)} & & \E\ar[r] & 
\ku\ar[r]\ar[d] & 0\\
& F\ar[d] & & & 0\\ & 0 }

\sepprop

\begin{subsub}\label{pr7}{\bf Lemme : }
Soient \ \m{\sigma=(\sigma_F,\sigma_G)} \ l'\'el\'ement de \
\m{\Ext^1_{\ko_X}(\ku,F)\oplus\Ext^1_{\ko_X}(\ku,G)} \ correspondant \`a la
ligne exacte du bas, et \ \m{\eta(\psi)\in\Ext^1_{\ko_X}(\ku,F)} celui
correspondant \`a la colonne exacte de droite. Alors on a
\[\eta(\psi)-\eta(0) \ = \ \ov{\psi}(\sigma_G) \ ,\]
\m{\ov{\psi}} d\'esignant le morphisme \ 
\m{\Ext^1_{\ko_X}(\ku,G)\to\Ext^1_{\ko_X}(\ku,F)} \ induit par $\psi$.
\end{subsub}

\begin{proof}
Soit
\xmat{\U_2\ar[r]^{f_2} & \U_1\ar[r]^{f_1} & \U_0\flon[r] & \ku}
une r\'esolution localement libre de $\ku$, telle que \
\m{\Ext^1_{\ko_X}(\U_0,G\oplus F)=\nsp}. De la suite exacte \ \m{0\to
\U_1/\imm(f_2)\to\U_0\to\ku\to 0}, on d\'eduit la suite exacte
\xmat{\Hom(\U_0,G\oplus F)\ar[r] & \Hom(\U_1/\imm(f_2),G\oplus
F)\ar[rr]^-{\delta=(\delta_G,\delta_F)} & & \Ext^1_{\ko_X}(\ku,G\oplus 
F)\ar[r] & 0 .}
Soient \ \m{\lambda:\U_1/\imm(f_2)\to G}, \m{\mu:\U_1/\imm(f_2)\to F} \ tels
que \ \m{\delta(\lambda,\mu)=\sigma}. Soient \Nligne 
\m{f'_1:\U_1/\imm(f_2)\to\U_0} \ le morphisme induit par \m{f_1} et
\[\gamma=(f'_1,\lambda,\mu):\U_1/\imm(f_2)\lra\U_0\oplus F\oplus G \ .\]
Alors on a \ \m{\E=\coker(\gamma)}, les morphismes \ \m{G\oplus F\to\E} \ et
\m{\E\to\ku} \ \'etant les morphismes \'evidents. Soit \
\m{\nu=\psi\circ\lambda+\mu\in\Hom(\U_1/\imm(f_2),F)}. Alors on a un
diagramme commutatif avec ligne exacte
\xmat{G\ar[rr]^-{(0,I_G,-\psi)} & & \U_0\oplus G\oplus F
\ar[rr]^-{(I,\psi+I_F)} & & \U_0\oplus F\\
& & & \U_1/\imm(f_2)\ar[ul]^{(f'_1,\lambda,\mu)}\ar[ur]_{(f'_1,\nu)}
}
d'o\`u on d\'eduit que \ \m{\ke_\psi=\coker(f'_1,\nu)}.
En faisant \m{\psi=0} on obtient \ \m{\eta(0)=\delta_F(\mu)}, et en g\'en\'eral
\[\eta(\psi) \ = \ \delta_F(\mu) + \delta_F(\psi\circ\lambda) \ = \
\eta(0)+ \ov{\psi}(\sigma_G) .\]
\end{proof}

\end{sub}

\sepsec

\section{Le faisceau canonique d'une courbe multiple primitive}\label{fai_can}

Soit \m{C_n} une courbe multiple primitive de multiplicit\'e $n$, de courbe
r\'eduite associ\'ee $C$ et de fibr\'e en droites sur $C$ associ\'e $L$. On
suppose comme dans \ref{prim_const} que \m{C_n} est obtenue en recollant les
vari\'et\'es \ \m{U_i\times Z_n} \ au moyen des automorphismes \m{\sigma_{ij}}
de \ \m{U_{ij}\times Z_n}. On note \m{\mathbf{U}_i} l'ouvert de \m{C_n}
correspondant \`a \m{U_i}.

Pour tout entier \m{k\geq 2}, on note \m{\ko_{0,k}} le faisceau structural de 
\m{C\times Z_k}.

Le faisceau canonique \m{\Omega_{C_n}} est {\em quasi localement libre} (cf.
\cite{dr4}, 3.4), localement isomorphe \`a \m{\ko_n\oplus\ko_{n-1}}. Soient
$P$ un point ferm\'e de $C$, \m{z\in\ko_{nP}} une \'equation de $C$ et \m{x\in
\ko_{nP}} au dessus d'un g\'en\'erateur de l'id\'eal maximal de \m{\ko_{C,P}}.
Alors \m{\Omega_{C_n,P}} est engendr\'e par \m{dx} (le facteur \m{\ko_n}) et
\m{dz} (le facteur \m{\ko_{n-1}})

\sepsub

\Ssect{Construction \`a partir de cocycles}{can_coc}

On va construire le faisceau canonique \m{\Omega_{C_n}} par la m\'ethode de
\ref{fai_rec}. On consid\`ere les isomorphismes \
\m{\Omega_{C_n\mid \mathbf{U}_i}\simeq\sigma_i^*(\Omega_{U_i\times Z_n})}. On
en d\'eduit les isomorphismes
\begin{equation}\label{equ1}
\Theta_{ij}:\Omega_{U_{ij}\times Z_n}\lra\sigma_{ij}^*(\Omega_{U_{ij}\times
Z_n}) .
\end{equation}
On a \ \m{\Omega_{U_{i}\times Z_n}\simeq\big(\ko_{0,n}\oplus \ko_{0,n-1}
\big)_{\mid U_{i}}}, avec des g\'en\'erateurs \m{dx}, \m{dt}, le premier
engendrant \m{\ko_{0,n\mid U_{i}}} et le second \m{\ko_{0,n-1\mid U_{i}}}.
Le faisceau \m{\Omega_{C_n}} est obtenu en recollant les faisceaux
\m{(\ko_{0,n}\oplus \ko_{0,n-1}\big)_{\mid U_{i}}} au moyen des isomorphismes
\m{\Theta_{ij}} par le proc\'ed\'e d\'ecrit dans \ref{fai_rec}.

Pour tout \ \m{f\in\ko_C(U_{ij})[t]/(t^n)}, on a
\[\Theta_{ij}(df) \ = \ d(f\circ\sigma_{ij}^{-1}) \ = \
d\big(\phi_{\mu_{ij},\nu_{ij}}(f)\big) .\]
Donc
\begin{eqnarray*}
\Theta_{ij}(dx) & = & d\big(\phi_{\mu_{ij},\nu_{ij}}(x)\big)\\
& = & d(x+\mu_{ij}t)\\
& = & (1+\frac{\partial\mu_{ij}}{\partial
x}t)dx+(\mu_{ij}+\frac{\partial\mu_{ij}}{\partial t}t)dt ,
\end{eqnarray*}
\begin{eqnarray*}
\Theta_{ij}(dt) & = & d\big(\phi_{\mu_{ij},\nu_{ij}}(t)\big)\\
& = & d(\nu_{ij}t)\\
& = & \frac{\partial\nu_{ij}}{\partial
x}t.dx+(\nu_{ij}+\frac{\partial\nu_{ij}}{\partial t}t)dt ,
\end{eqnarray*}
(exprim\'es dans le syst\`eme de g\'en\'erateurs \m{(dx,dt)} du
\m{\Omega_{U_{ij}\times Z_n}} de droite dans (\ref{equ1})). Mais, pour tenir
compte de \m{\sigma_{ij}^*} il faut appliquer \m{\phi_{\mu_{ij},\nu_{ij}}^{-1}}
aux coefficients de \m{dx}, \m{dt}. On obtient finalement que \m{\Theta_{ij}} 
est l'automorphisme de \ \m{\big(\ko_{0,n}\oplus \ko_{0,n-1}
\big)_{\mid U_{ij}}} d\'efini par la matrice
\begin{equation}\label{equ6}
\begin{pmatrix} 1+\phi'\big(\frac{\partial\mu_{ij}}{\partial x}
\big)\nu_{ji}t
& \phi'\big(\frac{\partial\nu_{ij}}{\partial x}\big)\nu_{ji}t \\
 & \\
-\frac{\mu_{ji}}{\nu_{ji}}+\phi'\big(\frac{\partial\mu_{ij}}{\partial t}\big)
\nu_{ji}t
& \frac{1}{\nu_{ji}}+\phi'\big(\frac{\partial\nu_{ij}}{\partial t}\big)
\nu_{ji}t
\end{pmatrix}
\end{equation}
en posant \ \m{\phi'=\phi_{\mu_{ij},\nu_{ij}}^{-1}} (on utilise
\ref{prop_phi}).

En particulier \m{\Omega_{C_n\mid C}} s'obtient en recollant les fibr\'es
\m{(\omega_C\oplus\ko_C)_{\mid U_i}} au moyen des automorphismes d\'efinis par
les matrices \m{\begin{pmatrix}1 & 0\\ \mu_{ij,0} & \nu_{ij,0}\end{pmatrix}}
(d'apr\`es \ref{prop_phi} on a \ \m{-\frac{\mu_{ji,0}}{\nu_{ji,0}}=
\mu_{ij,0}} \ et \ \m{\frac{1}{\nu_{ji,0}}=\nu_{ij,0}}).

\end{sub}

\sepsub

\Ssect{Faisceau canonique et prolongements de courbes multiples}{rest_can}

\begin{subsub}\label{pr1}{\bf Proposition :}
Le noyau du morphisme surjectif canonique \ \m{\rho_n:\Omega_{C_n\mid
C_{n-1}}\to \Omega_{C_{n-1}}} \ est isomorphe \`a \m{L^{n-1}}.
\end{subsub}

\begin{proof}
Soit $\ki$ le faisceau d'id\'eaux de \m{C_{n-1}} dans \m{C_n}. On a une suite
exacte
\xmat{\ki/\ki^2\ar[r]^-i & \Omega_{C_n\mid C_{n-1}}\ar[r] & 
\Omega_{C_{n-1}}\ar[r] & 0 .}
Soient $P$ un point ferm\'e de $C$ et \m{z\in\ko_{nP}} une \'equation de $C$.
Alors \ \m{\ki_{C,P}=(z^{n-1})} et \ \m{\imm(i_P)=(z^{n-2}dz)}. Dans les cartes
\m{(\ko_{0,n}\oplus \ko_{0,n-1}\big)_{\mid U_{i}}}, \m{\imm(i)} correspond aux
sous-faisceaux engendr\'es par \m{t^{n-2}dt}. D'apr\`es \ref{can_coc} ces
sous-faisceaux isomorphes \`a \m{\ko_{U_i}} se recollent par les isomorphismes
\[{\xymatrix@R=2pt{
\ko_{U_{ij}}\ar[r] & \ko_{U_{ij}} \\
\ u\fmaps[r] & \nu_{ij,0}^{n-1}u \quad ,}}\]
ce qui montre que \m{\imm(i)\simeq L^{n-1}}.
\end{proof}

\sepprop

On a donc une suite exacte
\begin{equation}\label{equ3}
0\lra L^{n-1}\lra\Omega_{C_n\mid C_{n-1}}\lra\Omega_{C_{n-1}}\lra 0 .
\end{equation}
Notons que \m{\Omega_{C_n\mid C_{n-1}}} est un fibr\'e vectoriel de rang 2 sur
\m{C_{n-1}}, et que \m{\rho_n} induit un isomorphisme \
\m{\Omega_{C_n\mid C}\simeq\Omega_{C_{n-1}\mid C}}.

R\'eciproquement, si $\E$ est un fibr\'e vectoriel de rang 2 sur \m{C_{n-1}},
et si \ \m{\phi:\E\to\Omega_{C_{n-1}}} \ est un morphisme surjectif, $\phi$
induit un isomorphisme entre les restrictions \`a $C$ et \m{\ker(\phi)\simeq
L^{n-1}} (\m{\ker(\phi)} se calcule en utilisant \ref{fib_det}).

On \'etudie maintenant les extensions \ \m{0\to 
L^{n-1}\to\ke\to\Omega_{C_{n-1}}\to 0}. D'abord localement :

\sepprop

\begin{subsub}\label{pr2}{\bf Lemme : } Soit $P$ un point ferm\'e de $C$.

1 - On a \ \m{\Ext^1_{\ko_{n-1,P}}(\ko_{n-1,P}\oplus\ko_{n-2,P},\ko_{C,P})
\simeq\ko_{C,P}} .

2 - Soit \ \m{0\to\ko_{C,P}\to M\to\ko_{n-1,P}\oplus\ko_{n-2,P}\to 0} \ une
suite exacte de \m{\ko_{n-1,P}}-modules, associ\'ee \`a \m{\alpha\in\ko_C}.
Alors on a \ \m{M\simeq 2\ko_{n-1,P}} \ si et seulement si $\alpha$ est
inversible.
\end{subsub}
On a \ \m{\Ext^1_{\ko_{n-1,P}}(\ko_{n-1,P},\ko_{C,P})=\nsp}, donc il suffit de
prouver :
\begin{enumerate}
\item[(i)] On a \ $\Ext^1_{\ko_{n-1,P}}(\ko_{n-2,P},\ko_{C,P})\simeq\ko_{C,P}$.
\item[(ii)] Soit \ \m{0\to\ko_{C,P}\to N\to\ko_{n-2,P}\to 0} \ une suite exacte
de \m{\ko_{n-1,P}}-modules, associ\'ee \`a \m{\alpha\in\ko_C}. Alors on a \
\m{N\simeq\ko_{n-1,P}} \ si et seulement si $\alpha$ est inversible.
\end{enumerate}
\begin{proof}
Soit \m{z\in\ko_{n-1,P}} un g\'en\'erateur de l'id\'eal de $C$. Alors (i)
d\'ecoule imm\'ediatement de la r\'esolution libre de \m{\ko_{n-2,P}}
\xmat{\ldots\ko_{n-1,P}\ar[r]^{z} & \ko_{n-1,P}\ar[r]^{z^{n-2}} &
\ko_{n-1,P}\ar[r] & \ko_{n-2,P} .}
D\'emontrons maintenant (ii). On a une suite exacte
\xmat{\ko_{C,P}\ar[r]^-\delta & \Ext^1_{\ko_{n-1,P}}(\ko_{n-2,P},\ko_{C,P})
=\ko_{C,P}\ar[r] & \Ext^1_{\ko_{n-1,P}}(N,\ko_{C,P}) ,}
obtenue en appliquant \m{\Hom(-,\ko_{C,P})} \`a la suite exacte de (ii). Le
morphisme $\delta$ est la multiplication par $\alpha$.

Si $\alpha$ n'est pas inversible, \m{\coker(\delta)\not=0}, donc \
\m{\Ext^1_{\ko_{n-1,P}}(N,\ko_{C,P})\not=0}, et \m{N\not=\ko_{n-1,P}}.

La suite exacte canonique
\xmat{0\ar[r] & \ko_{C,P}\simeq(z^{n-2})\ar[r] & \ko_{n-1,P}\ar[r]^1 & 
\ko_{n-2,P}\ar[r] & 0}
est donc associ\'ee \`a un \'el\'ement inversible de \m{\ko_{C,P}}.

Supposons que $\alpha$ est inversible. La suite exacte
\xmat{0\ar[r] & \ko_{C,P}\simeq(z^{n-2})\ar[r] & 
\ko_{n-1,P}\ar[r]^{\alpha_0/\alpha} & \ko_{n-2,P}\ar[r] & 0}
est associ\'ee \`a $\alpha$, donc \m{N\simeq\ko_{n-1,P}}.
\end{proof}

\sepprop

On a une suite exacte canonique
\xmat{0\ar[r] & H^1(\HHom(\Omega_{C_{n-1}},L^{n-1}))\ar[r] & \qquad\qquad\qquad
\qquad\qquad\qquad\qquad\qquad}
\xmat{\ar[r] & \Ext^1_{\ko_{n-1}}(\Omega_{C_{n-1}},L^{n-1})\ar[r]^-\pi &
H^0(\EExt^1(\Omega_{C_{n-1}},L^{n-1}))\ar[r] & 0 .}
d'apr\`es la suite spectrale des Ext (cf. \cite{go}, 7.3).

\sepprop

\begin{subsub}\label{pr3}{\bf Lemme : } On a \
\m{\EExt^1(\Omega_{C_{n-1}},L^{n-1})\simeq\ko_C} \ et \
\m{\HHom(\Omega_{C_{n-1}},L^{n-1})\simeq E\ot L^{n-1}} .
\end{subsub}

\begin{proof} Soit $\L$ un fibr\'e en droites sur \m{C_{n-1}} tel que
\m{\L_{\mid C}\simeq L} (cf. \cite{dr2}, 3.1.1). On a une r\'esolution
localement libre de $L$
\[\cdots\lra\L^n\lra\L^{n-1}\lra L^{n-1} ,\]
d'o\`u, en utilisant la suite exacte (\ref{equ3}) la r\'esolution localement 
libre de \m{\Omega_{C_{n-1}}}
\xmat{\cdots\ar[r] & \L^n\ar[r]^-{f_1} & \L^{n-1}\ar[r]^-{f_0} &
\Omega_{C_n\mid C_{n-1}}\ar[r] & \Omega_{C_{n-1}}}
avec laquelle on peut calculer \m{\EExt^1(\Omega_{C_{n-1}},L^{n-1})}. Le
premier r\'esultat d\'ecoule du fait que \m{f_{0\mid C}} et \m{f_{1\mid C}}
s'annulent.

Le second est imm\'ediat car \ \m{\Omega_{C_{n-1}\mid C}=\Omega_{C_2\mid
C}=E^*}.
\end{proof}

\sepprop

Soit \m{\sigma(C_n)} l'\'el\'ement de 
\m{\Ext^1_{\ko_{n-1}}(\Omega_{C_{n-1}},L^{n-1})} correspondant \`a la suite
exacte (\ref{equ3}). On a \ \m{H^0(\EExt^1(\Omega_{C_{n-1}},L^{n-1}))=\C} \
d'apr\`es le lemme \ref{pr3}, et d'apr\`es le lemme \ref{pr2} on peut supposer
que \ \m{\pi(\sigma(C_n))=1} .

Soit maintenant \m{C'_n} une autre courbe multiple primitive de multiplicit\'e
$n$ extension de \m{C_{n-1}}. On a alors \ \m{\sigma(C'_n)-\sigma(C_n)\in
H^1(E\ot L^{n-1})}.

On conserve les notations de \ref{prolong}.

D'apr\`es \ref{can_coc}, \m{\Omega_{C_n\mid C_{n-1}}} est obtenu en recollant
les faisceaux \m{\big(\ko_{0,n-1}\oplus \ko_{0,n-1}\big)_{\mid U_{ij}}} au
moyen des automorphismes d\'efinis par les matrices (\ref{equ6}), qu'on note
\m{M_{ij}}. On note \m{M'_{ij}} les matrices (\ref{equ6}) pour
\m{\Omega_{C'_n\mid C_{n-1}}}. Un calcul simple montre que
\[M'_{ij}-M_{ij} \ = \ \begin{pmatrix} 0 & 0 \\
(n-1)\alpha_{ij}\nu_{ij,0}^{2-n}t^{n-2} & 
(n-1)\beta_{ij}\nu_{ij,0}^{2-n}t^{n-2}
\end{pmatrix} \]
On en d\'eduit ais\'ement, avec la discussion de \ref{prolong} la

\sepprop

\begin{subsub}\label{pr6}{\bf Proposition : }
Pour tout \ \m{u\in H^1(E\ot L^{n-1})} \ on a \ \m{\sigma(C_n(u))-\sigma(C_n)
=(n-1)u} .
\end{subsub}

\end{sub}

\sepsec

\section{Courbes multiples et familles de courbes lisses}\label{demo}

\Ssect{Courbes multiples provenant d'une famille de courbes lisses}{cm_fam}

Soit \m{\pi:\kd\to S} une famille de plate de courbes projectives 
lisses
irr\'eductibles param\'etr\'ee par une vari\'et\'e lisse irr\'eductible $S$ de 
dimension \m{d>0}. Soient \m{s_0} un point ferm\'e de $S$ et 
\m{C=\pi^{-1}(s_0)}.

Soit \m{n\geq 3} un entier. On pose \ \m{Z_n=\spec({\C[t]/(t^n)})}. Soit \ 
\m{\phi:Z_n\hookrightarrow S} \ un plongement tel que l'image du point ferm\'e 
de \m{Z_n} soit \m{s_0}. Alors \ \m{C_n=\pi^{-1}(Z_n)} \ est une courbe 
multiple primitive de multiplicit\'e $n$, de courbe lisse associ\'ee $C$, et 
le faisceau d'id\'eaux de $C$ dans \m{C_n} est \m{\ko_{n-1}}.
Soient \m{t_1,\ldots,t_d} des \'el\'ements de l'id\'eal maximal 
\m{\mathbf{m}_{s_0}} de \m{\ko_{S,s_0}} formant une base de 
\m{\mathbf{m}_{s_0}/\mathbf{m}_{s_0}^2}. Soit \ 
\m{\Phi:\ko_{S,s_0}\to\C[t]/(t^n)} \ le morphisme correspondant \`a $\phi$. Il 
est enti\`erement d\'etermin\'e par \m{\Phi(t_1)\ldots,\Phi(t_d)}. En 
changeant \m{t_1,\ldots,t_d} on se ram\`ene ais\'ement au cas o\`u \ 
\m{\Phi(t_2)=\cdots\Phi(t_d)=0} \ et \ \m{\Phi(t_1)=t}. Soit \m{\ki_{C_n}} le
faisceau d'id\'eaux de \m{C_n} dans $\kd$. On a \ \m{\ki_{C_n}=\langle t_1^n,
t_2,\cdots,t_d\rangle}.

Soit \ \m{\phi':Z_n\hookrightarrow S} \ un autre plongement tel que \ 
\m{\Phi'_{\mid Z_{n-1}}=\Phi_{\mid Z_{n-1}}}. Si \ 
\m{\Phi':\ko_{S,s_0}\to\C[t]/(t^n)} \ est le morphisme correspondant, il 
existe \m{\alpha_1,\alpha_2,\ldots,\alpha_d\in\C} tels que
\[\Phi'(t_1)=t+\alpha_1t^{n-1},\quad\Phi'(t_2)=\alpha_2t^{n-1},\ldots,
\Phi'(t_d)=\alpha_dt^{n-1} .\]
Si \m{C'_n} est la courbe multiple d\'efinie par \m{\phi'}, on a
\[\ki_{C'_n} \ = \ \langle t_1^n,t_2-\alpha_2t_1^{n-1},\ldots,
t_d-\alpha_dt_1^{n-1}\rangle \ . \]

\end{sub}

\sepsub

\Ssect{Faisceaux canoniques}{f_cano}

On a une suite exacte canonique
\xmat{\ki_{C_n}/\ki_{C_n}^2\ar[r]^-{A_n} & \Omega_{\kd\mid C_n}\ar[r] & 
\Omega_{C_n}\ar[r] & 0 .}
Le faisceau \m{\Omega_{\kd\mid C_n}} est localement libre de rang \m{d+1} et 
la structure de \m{\Omega_{C_n}} est donn\'ee dans \ref{fai_can}.

On pose \ \m{\Gamma=\langle t_2,\ldots,t_d\rangle\subset 
\mathbf{m}_{s_0}/\mathbf{m}_{s_0}^2}. Soient \m{P\in C} et \m{z\in\ko_{\kd P}} 
au dessus d'un g\'en\'erateur de l'id\'eal maximal de \m{\ko_{C,P}}. Alors 
\m{(dz, dt_1,\ldots,dt_d)} est une base de \m{\Omega_{\kd\mid C_n,P}}. On a 
donc
\[\imm(A_n) \ = \ \langle t_1^{n-1}dt_1,dt_2,\ldots,dt_d\rangle \
\simeq \ \ko_C\oplus(\ko_n\ot\Gamma) \ . \]
On en d\'eduit la suite exacte
\xmat{0\ar[r] & \ko_{n-1}\ot\Gamma\ar[r]^-\beta & \Omega_{\kd\mid 
C_{n-1}}\ar[r] & \Omega_{C_n\mid C_{n-1}}\ar[r] & 0 \ .}
On a \ \m{\imm(\beta)= \langle dt_2,\ldots,dt_d\rangle} . On obtient le 
diagramme commutatif avec lignes et colonnes exactes
\xmat{& & & & 0\ar[d] \\ & 0\ar[d] & & & \ko_C\ar[d] \\
0\ar[r] & G\ar[rr]^-{\beta(I_G,0)}\ar[d]^{(I_G,0)} & & 
\Omega_{\kd\mid C_{n-1}}\fleq[d]\ar[r] & \Omega_{C_n\mid C_{n-1}}\ar[r]\ar[d]& 
0\\
0\ar[r] & G\oplus\ko_C\ar[rr]^-\beta\ar[d]^{(0,I_{\ko_C})} & & \Omega_{\kd\mid 
C_{n-1}}\ar[r] & \Omega_{C_{n-1}}\ar[r]\ar[d] & 0\\
& \ko_C\ar[d] & & & 0\\ & 0 }
avec \ \m{G=\ko_{n-1}\ot\Gamma}, \m{\imm(\beta)=\langle 
t_1^{n-2}dt_1,dt_2,\ldots,dt_d\rangle} .
Pour la courbe \m{C'_n} on a le diagramme suivant
\xmat{& & & & 0\ar[d] \\ & 0\ar[d] & & & \ko_C\ar[d] \\
0\ar[r] & G\ar[rr]^-{\beta(I_G,-\psi)}\ar[d]^{(I_G,-\psi)} & & 
\Omega_{\kd\mid C_{n-1}}\fleq[d]\ar[r] & \Omega_{C'_n\mid 
C_{n-1}}\ar[r]\ar[d]& 0\\
0\ar[r] & G\oplus\ko_C\ar[rr]^-\beta\ar[d]^{(\psi,I_{\ko_C})} & & 
\Omega_{\kd\mid 
C_{n-1}}\ar[r] & \Omega_{C_{n-1}}\ar[r]\ar[d] & 0\\
& \ko_C\ar[d] & & & 0\\ & 0 }
o\`u \ 
\m{\psi=((n-1)\alpha_2,\ldots,(n-1)\alpha_d):\ko_{n-1}\ot\Gamma\to\ko_C}.
\end{sub}

\sepsub

\Ssect{D\'emonstration du th\'eor\`eme \ref{theo}}{dem_theo}

On utilise le r\'esultat suivant :

\sepprop

\begin{subsub}\label{pr8}{\bf Proposition : }
Il existe une famille de courbes lisses \ \m{\pi:\kd\to S} \ param\'etr\'ee 
par une vari\'et\'e lisse \m{S}, telle qu'il existe \m{s_0\in S} tel que 
\m{\kd_{s_0}\simeq C}, et que le morphisme de Koda\"\i ra-Spencer \ 
\m{KS:T_{s_0}S\to H^1(T_C)} \ soit surjectif.
\end{subsub}

\begin{proof} Soit $\kl$ un fibr\'e en droites tr\`es ample sur $C$, et \ 
\m{C\hookrightarrow\P_n=\P(H^0(\kl)^*)} \ le plongement induit. On suppose 
que \ \m{H^1(\kl)=\nsp}. En consid\'erant la suite exacte canonique
\[0\lra \ko_{\P_n}\lra\ko_{\P_n}(1)\ot H^0(\ko_{\P_n}(1))^*\lra T\P_n
\lra 0 \]
on voit que \ \m{H^1(T{\P_n}_{\mid C})=\nsp}. Soient $P$ le polyn\^ome de 
Hilbert de \m{\ko_C} et \m{\Hilb^P(\P_n)} le sch\'ema de Hilbert correspondant.
D'apr\`es les propri\'et\'es diff\'erentielles de ce sch\'ema (cf. 
\cite{grot}), \m{\Hilb^P(\P_n)} est lisse au point \m{s_0} correspondant \`a 
$C$, et le sch\'ema universel $\kd$ au voisinage de \m{s_0} est la famille de 
courbes lisses recherch\'ee.
\end{proof}

\sepprop

On suppose donn\'ee une famille de courbes lisses $\kd$ comme dans la 
proposition \ref{pr8}. Soit $X$ une sous-vari\'et\'e lisse de $S$ contenant 
\m{s_0}, et \m{\kd_X} l'image inverse de $X$ dans $\kd$. Alors on a une suite 
exacte canonique
\begin{equation}\label{equ7}0\lra\ko_C\ot[T_{s_0}X]^*\lra\Omega_{X\mid 
C}\lra\omega_C\lra 0 ,\end{equation}
et l'\'el\'ement associ\'e de \ 
\m{\Ext^1_{\ko_C}(\omega_C,\ko_C\ot[T_{s_0}X]^*)=\Hom(T_{s_0}X,H^1(T_C))} \ 
n'est autre que le morphisme de Koda\"\i ra-Spencer de \m{\kd_X} en \m{s_0}, 
qui est la restriction \`a \m{T_{s_0}X} de celui de $\kd$.

Le th\'eor\`eme \ref{theo} se d\'emontre par r\'ecurrence sur 
$n$. 

D'apr\`es \ref{cm_fam}, il suffit de montrer que pour toute courbe multiple 
primitive \m{C_n} de multiplicit\'e $n$ telle que le faisceau d'id\'eaux de 
$C$ dans \m{C_n} soit le fibr\'e trivial sur \m{C_{n-1}}, il existe un 
plongement \ \m{\phi:Z_n\hookrightarrow S} \ tel que l'image du point ferm\'e 
de \m{Z_n} soit \m{s_0} et que \m{C_n\simeq\pi^{-1}(Z_n)}.

Traitons d'abord le cas \m{n=2}. On suppose que \m{C_2} est une courbe double 
de courbe lisse associ\'ee $C$, telle que \m{\ki_C=L} soit trivial sur $C$.
Si \m{C_2} est triviale, la famille triviale de courbes \m{C\times\C} r\'epond
aux conditions du th\'eor\`eme \ref{theo}. Supposons que \m{C_2} n'est pas 
triviale. Soit \m{D\subset H^1(T_C)} la droite engendr\'ee par l'\'el\'ement 
de \m{H^1(T)} correspondant \`a l'extension
\begin{equation}\label{equ8}0\lra\ko_C\lra\Omega_{C_2\mid C}\lra\omega_C\lra 0 
.\end{equation}
Soit \m{Y\subset S} une courbe lisse passant par \m{s_0} et telle que \ 
\m{KS(T_{s_0}Y)=D}. Soit \m{C'_2} le second voisinage infinit\'esimal de $C$ 
dans $\kd_Y$. Si \m{z\in\ko_{Y,s_0}} est un g\'en\'erateur de l'id\'eal 
maximal, le faisceau d'id\'eaux de \m{C'_2} dans \m{\kd_Y} est engendr\'e par 
\m{z^2}. On a donc \m{\Omega_{C'_2\mid C}=\Omega_{\kd_Y\mid C}}, et la suite 
exacte (\ref{equ7}) pour $Y$ est la m\^eme que la suite exacte (\ref{equ8}) 
pour \m{C'_2}. Il en d\'ecoule que la droite de \m{H^1(T)} d\'efinie par 
l'extension
\[0\lra\ko_C\lra\Omega_{C'_2\mid C}\lra\omega_C\lra 0\]
est \'egale \`a $D$. D'apr\`es \ref{mult_2}, \m{C'_2} est isomorphe \`a 
\m{C_2}, ce qui d\'emontre le th\'eor\`eme \ref{theo} pour \m{n=2}.

\medskip

Supposons le th\'eor\`eme \ref{theo} d\'emontr\'e pour les courbes de 
multiplicit\'e \m{n-1}. Soit \m{C_{n-1}} une courbe multiple 
primitive de multiplicit\'e \m{n-1} telle que le faisceau d'id\'eaux de 
$C$ dans \m{C_{n-1}} soit le fibr\'e trivial sur \m{C_{n-2}}. Il faut montrer 
que toute extension \m{C_n^0} de \m{C_{n-1}} en courbe multiple primitive de 
multiplicit\'e $n$ telle que le faisceau d'id\'eaux de $C$ dans \m{C_n^0} soit 
le fibr\'e trivial sur \m{C_{n-1}} provient d'une famille de courbes lisses.
On peut supposer que \m{C_{n-1}} provient d'un plongement \ 
\m{Z_{n-1}\hookrightarrow S} \ tel que l'image du point ferm\'e 
de \m{Z_n} soit \m{s_0}, et que ce plongement est \'etendu \`a un plongement  
\ \m{\phi:Z_n\hookrightarrow S} \ correspondant \`a la courbe multiple 
\m{C_n}, extension de \m{C_{n-1}}. On reprend les notations de \ref{cm_fam} et 
\ref{f_cano}.

On a un diagramme commutatif avec lignes et colonnes exactes
\xmat{0\ar[r] & G\oplus\ko_C\ar[r]^-\alpha\ar[d]^A & \Omega_{\kd\mid C_{n-1}}
\ar[r]\ar[d] & \Omega_{C_{n-1}}\ar[r]\ar[d]^r & 0\\
0\ar[r] & \ko_C\ot[T_{s_0}S]^*\ar[r] & \Omega_{\kd\mid C}\ar[r] & 
\omega_C\ar[r] 
& 0}
Le morphisme $A$ est nul sur \m{\ko_C}, et provient de l'inclusion 
\m{\Gamma\subset [T_{s_0}S]^*} sur $G$.

La suite exacte du haut est associ\'ee \`a
\[(\sigma_G,\sigma_{\ko_C}) \ \in \ \Ext^1_{\ko_{n-1}}(\Omega_{\kd\mid 
C_{n-1}},G\oplus\ko_C) . \]
Il existe des sections locales de \m{\alpha_{\mid G}}, donc on a \ 
\m{\sigma_G\in H^1(\HHom(\Omega_{\kd\mid C_{n-1}},G))}.

La suite exacte du bas est associ\'ee \`a
\[\sigma' \ \in \ \Ext^1_{\ko_{n-1}}(\omega_C,\ko_C\ot[T_{s_0}S]^*) .\]
Puisque cette suite est localement scind\'ee, on a
\[\sigma' \ \in \ H^1(\HHom(\omega_C,\ko_C\ot[T_{s_0}S]^*)) \ = \
\Hom(T_{s_0}S,H^1(T_C)) ,\]
et \m{\sigma'} n'est autre que le morphisme de Koda\"\i ra-Spencer
\m{KS:T_{s_0}S\to H^1(T)} \ de $\kd$ en \m{s_0}.

De \m{\sigma'} on d\'eduit
\[r(\sigma') \ \in \ H^1(\HHom(\Omega_{\kd\mid C_{n-1}},\ko_C\ot[T_{s_0}S]^*))
\ = \ \Hom(T_{s_0}S, H^1(E)) ,\]
qui est le compos\'e
\xmat{T_{s_0}S\ar[r]^-{KS} & H^1(T_C)\flon[r] & H^1(T_C)/D\flinc[r] & H^1(E)}
(cf. \ref{ic_tr}).

De \m{\sigma_G} on d\'eduit
\[A(\sigma_G) \ \in \ H^1(\HHom(\omega_C,\ko_C\ot[T_{s_0}S]^*)) \ = \
\Hom(T_{s_0}S,H^1(T_C)) .\]
D'apr\`es le diagramme commutatif pr\'ec\'edent, on a \ 
\m{r(\sigma')=A(\sigma_G)}.
En particulier, \m{A(\sigma_G)} se factorise par \m{\Gamma^*} :
\xmat{T_{s_0}S\ar[r]^-{KS}\flon[dr] & H^1(T_C)\flon[r] & H^1(T_C)/D\\
& \Gamma^*\ar[ur]^\lambda}
et $\lambda$ est surjective.

On consid\`ere maintenant la courbe \m{C'_n} de \ref{f_cano}. Elle est de la 
forme \ \m{C'_n=C_n(u)}, avec \ \m{u\in H^1(T_C)/D\subset H^1(E)} (cf. 
\ref{rest_can}). D'apr\`es le lemme \ref{pr7} et la proposition 3.2.4, on a 
\m{u=\frac{\lambda(\psi)}{n-1}}.

D'apr\`es 2.4.5, il existe \m{u_0\in H^1(T_C)/D} tel que
\m{C_n^0=C_n(u_0)}. Puisque $\lambda$ est surjective, on peut choisir $\psi$
tel que \m{\frac{\lambda(\psi)}{n-1}=u_0}, ce qui prouve que \m{C_n^0} 
provient bien d'une famille de courbes lisses.

\end{sub}

\end{document}